\documentclass[12pt]{article}
\usepackage{amsmath,amsthm,
amssymb,extsizes}
\title{Tridiagonal
canonical matrices of
bilinear or
sesquilinear forms and
of pairs of symmetric,
skew-symmetric, or
Hermitian forms%
\thanks{Preprint RT-MAT 2006-16, Universidade de Sao Paulo, 2006, 21 p.}}

\author{Vyacheslav Futorny%
\thanks{Partially
supported by CNPq,
processo
307812/2004-9, and by
FAPESP, processo
2005/60337-2.}
\\Department of Mathematics, University of
S\~{a}o Paulo\\
S\~{a}o Paulo, Brazil,
futorny@ime.usp.br
 \and
Roger A. Horn%
\\Department of Mathematics, University of
Utah\\ Salt Lake City,
Utah 84112-0090,
rhorn@math.utah.edu
\and
Vladimir V. Sergeichuk%
\thanks{Corresponding
author. Partially
supported by FAPESP,
processo 05/59407-6.
This author is
grateful to the
University of S\~ao
Paulo for the
hospitality.}
\\ Institute of
Mathematics,
Tereshchenkivska 3\\
Kiev, Ukraine,
sergeich@imath.kiev.ua}
\date{}

\DeclareMathOperator{\rank}{rank}
\DeclareMathOperator{\diag}{diag}

\begin{document}
\maketitle

\renewcommand{\le}{\leqslant}
\renewcommand{\ge}{\geqslant}

\newtheorem{theorem}{Theorem}[section]
\newtheorem{lemma}{Lemma}[section]

\theoremstyle{remark}
\newtheorem{remark}{Remark}[section]

\newcommand{\dia}{\,\diagdown\,}

\newcommand{\ddd}{
\text{\begin{picture}(12,8)
\put(-2,-4){$\cdot$}
\put(3,0){$\cdot$}
\put(8,4){$\cdot$}
\end{picture}}}

\begin{abstract}
Tridiagonal canonical
forms of square
matrices under
congruence or
*congruence, pairs of
symmetric or
skew-symmetric
matrices under
congruence, and pairs
of Hermitian matrices
under *congruence are
given over an
algebraically closed
field of
characteristic
different from~$2$.


\end{abstract}

\section{Introduction}

We give tridiagonal
canonical forms of
matrices of
\begin{itemize}
  \item[\rm(i)]
bilinear forms and
sesquilinear forms,

  \item[\rm(ii)]
pairs of forms, in
which each form is
either symmetric or
skew-symmetric, and

  \item[\rm(iii)]
pairs of Hermitian
forms
\end{itemize}
over an algebraically
closed field of
characteristic
different from $2$.
Our canonical forms
are direct sums of
matrices or pairs of
matrices of the form
\begin{equation}
\label{kdu}
\begin{bmatrix}
\varepsilon &a&&&&0\\
a'&0&b\\
&b'&0&a\\
&&a'&0&b\\
&&&b'&0&\ddots\\
0&&&&\ddots&\ddots
  \end{bmatrix};
\end{equation}
they employ relatively
few different types of
canonical direct
summands.

Let $\mathbb F$ be a
field of
characteristic
different from $2$.
The problem of
classifying bilinear
or sesquilinear forms
over $\mathbb F$ was
reduced by Gabriel,
Riehm, and
Shrader-Frechette
\cite{gab_form,rie,rie1}
to the problem of
classifying Hermitian
forms over finite
extensions of $\mathbb
F$. In \cite{ser_izv}
this reduction was
extended to
selfadjoint
representations of
linear categories with
involution, and
canonical matrices of
(i)--(iii) were
obtained over $\mathbb
F$ up to
classification of
Hermitian forms over
finite extensions of
$\mathbb F$. Canonical
matrices were found in
a simpler form in
\cite{hor-ser2} when
$\mathbb F=\mathbb C$.
Canonical matrices of
bilinear forms over an
algebraically closed
field of
characteristic $2$
were given in
\cite{ser_char2}. The
problem of classifying
pairs of symmetric,
skew-symmetric, or
Hermitian forms was
studied by many
authors; we refer the
reader to Thompson's
classical work
\cite{thom} with a
bibliography of 225
items, and to the
recent papers by
Lancaster and Rodman
\cite{lan-rod,lan-rod1}.

Each $n\times n$
matrix $A$ over
$\mathbb F$ defines a
bilinear form $x^TAy$
on ${\mathbb F}^n$. If
$\mathbb F$ is a field
with a fixed
nonidentity involution
$a\mapsto \bar a$,
then $A$ defines a
sesquilinear form
$\bar x^TAy$ on
${\mathbb F}^n$. Two
square matrices $A$
and $A'$ give the same
bilinear
(sesquilinear) form
with respect to
different bases if and
only if they are
\emph{congruent}
({\it\!*congruent});
this means that there
is a nonsingular $S$
such that $S^TAS=A'$
($S^*AS=A'$ with
$S^*:=\bar S^T$,
respectively). Two
matrix pairs $(A,B)$
and $(A',B')$ are
\emph{congruent}
(\!\emph{*congruent})
if there is a
nonsingular $S$ such
that $S^TAS=A'$ and
$S^TBS=B'$ ($S^*AS=A'$
and $S^*BS=B'$,
respectively). A
matrix $A$ is
\emph{Hermitian} if
$A=A^*$.

Thus, the canonical
form problem for
(i)--(iii) is the
canonical form problem
for
\begin{itemize}
  \item[\rm(i$'$)]
matrices under
congruence or
*congruence (their
tridiagonal canonical
matrices are given in
Theorems \ref{t_matr}
and \ref{t_matr1});

  \item[\rm(ii$'$)]
pairs of matrices
under congruence, in
which each matrix is
either symmetric or
skew-symmetric
(Theorems
\ref{th1a}--\ref{t_skew});
and

  \item[\rm(iii$'$)]
pairs of Hermitian
matrices under
*congruence (Theorem
\ref{t_her}).
\end{itemize}

The problem of finding
tridiagonal canonical
forms of (ii$'$) or
(iii$'$) is connected
with the problem of
tridiagonalizing
matrices by orthogonal
or unitary similarity:
two pairs $(I_n,B)$
and $(I_n,B')$ are
congruent or
*congruent if and only
if $B$ and $B'$ are
orthogonally or
unitarily similar,
respectively. The
well-known algorithm
for reducing symmetric
real matrices to
tridiagonal form by
orthogonal similarity
\cite[Section 5]{wil}
can not be extended to
symmetric complex
matrices. However,
Ikramov
\cite{ikr_trid} showed
that every symmetric
complex matrix is
orthogonally similar
to a tridiagonal
matrix. Each $4\times
4$ complex matrix is
unitarily similar to a
tridiagonal matrix
\cite{dav-dok,pati},
but there is a
$5\times 5$ matrix
that is not unitarily
similar to a
tridiagonal matrix
\cite{dav-dok,fon-wu,lon,stu}.

Our paper was inspired
by
\cite{dok-zhao_tridiag},
in which
\raisebox{1.5pt}{-}\!\!Dokovi\'{c}
and Zhao gave a
tridiagonal canonical
form of symmetric
matrices for
orthogonal similarity
over an algebraically
closed field of
characteristic
different from $2$ (we
use it in Theorem
\ref{th1b} of our
paper). In a
subsequent article,
and for the same type
of field,
\raisebox{1.5pt}{-}\!\!Dokovi\'{c},
Rietsch, and Zhao
\cite{dok-zhao_skew}
found a $4$-diagonal
canonical form of
skew-symmetric
matrices for
orthogonal similarity.

Matrix pairs $(A,B)$
and $(A',B')$ are
\emph{equivalent} if
there are nonsingular
$R$ and $S$ such that
$RAS=A'$ and $RBS=B'$.
We denote equivalence
of pairs by $\approx$.
Kronecker's theorem on
pencils of matrices
\cite[Section XII,
Theorem 5]{gan}
ensures that each pair
of matrices of the
same size is
equivalent to a direct
sum, determined
uniquely up to
permutation of
summands, of pairs of
the form
\begin{equation*}\label{ksca}
(I_n,J_n(\lambda)),\quad
(J_{n}(0),I_{n}),\quad
(F_n,G_n),\quad
(F_n^T,G_n^T),
\end{equation*}
in which $I_n$ is the
$n\times n$ identity
matrix,
\[
J_n(\lambda) :=
\begin{bmatrix}
  \lambda&1&&0\\
  &\lambda&\ddots&\\
  &&\ddots&1\\
0&&&\lambda
\end{bmatrix}
\text{ is $n$-by-$n$,}
\]
and
\begin{equation*}\label{1.4}
F_n:=\begin{bmatrix}
1&0&&0\\&\ddots&\ddots&\\0&&1&0
\end{bmatrix}\text{
and }\
G_n:=\begin{bmatrix}
0&1&&0\\&\ddots&\ddots&\\0&&0&1
\end{bmatrix}\text{
are $n$-by-$(n+1)$.}
\end{equation*}
Thus, $F_0=G_0$ is the
$0$-by-$1$ matrix,
which represents the
linear mapping
${\mathbb F}\to 0$.

In the following two
theorems (proved in
Sections \ref{s_pr}
and \ref{s_pr1}) we
give tridiagonal
canonical forms of a
square matrix $A$
under congruence and
*congruence. We also
give the Kronecker
canonical form of
$(B^T,B)$ and,
respectively,
$(B^*,B)$ for each
canonical direct
summand $B$, which
permits us to
construct the
canonical form of $A$
for congruence using
the Kronecker
canonical form of
$(A^T,A)$, and to
construct, up to signs
of the direct
summands, the
canonical form of $A$
for *congruence using
the Kronecker
canonical form of
$(A^*,A)$.

\begin{theorem}
\label{t_matr}

{\rm(a)} Each square
matrix $A$ over an
algebraically closed
field $\mathbb F$ of
characteristic
different from $2$ is
congruent to a direct
sum, determined
uniquely up to
permutation of
summands, of
tridiagonal matrices
of three types:
\begin{equation}\label{cm1}
\begin{bmatrix}
    0&1 &&&0\\
     \lambda&0&1 \\
    &\lambda &0&\ddots&\\
    &&\ddots&\ddots&1 \\
    0&&&\lambda &0
  \end{bmatrix}_{2k},\qquad
  \lambda
\in\mathbb F,\
\lambda\ne \pm 1,
\end{equation}
in which each nonzero
$\lambda$ is
determined up to
replacement by
$\lambda^{-1}$
$($i.e., the matrices
\eqref{cm1} with
$\lambda$ and
$\lambda^{-1}$ are
congruent$)$;
\begin{equation}\label{cm2}
  \begin{bmatrix}
    \varepsilon &1&&&&0\\
    -1&0&1&\\
    &1&0&1&\\
  &&-1&0&1&\\
&&&1&0&\ddots&\\
  0&&&&  \ddots&\ddots
  \end{bmatrix}_{n},
    \quad
  \begin{aligned}
&\text{$\varepsilon
=1$ if $n$ is a
multiple
of $4$,}\\
&\text{$\varepsilon
\in\{0,1\}$
otherwise};
  \end{aligned}
\end{equation}
and
\begin{equation}\label{cm3}
\begin{bmatrix}
    0 &1&&&&0\\
    1&0&1&\\
    &-1&0&1&\\
    &&1&0&1\\
&&&-1&0&\ddots&\\
  0&&&&  \ddots&\ddots\\
  \end{bmatrix}_{4k}.
  \qquad\qquad
    \qquad\qquad
      \qquad
\end{equation}
The subscripts $2k,\
n$, and $4k$ $($with
$k,n\in\mathbb N)$
designate the sizes of
the corresponding
matrices.

{\rm(b)} The direct
sum asserted in
{\rm(a)} is determined
uniquely up to
permutation of
summands by the
Kronecker canonical
form of $(A^T,A)$ for
equivalence. For each
direct summand $B$ of
types
\eqref{cm1}--\eqref{cm3},
the Kronecker
canonical form of
$(B^T,B)$ is given in
the following table:
\begin{equation}\label{tab1}
\renewcommand{\arraystretch}{1.2}
\begin{tabular}{|c|c|}
\hline $B$& Kronecker
canonical form of $
(B^T,B)$\\
\hline\hline Matrix
\eqref{cm1}&
$(I_k,J_k(\lambda))\oplus
(J_k(\lambda),I_k)$\\\hline
Matrix \eqref{cm2}&
$(F_k,G_k)\oplus
(F_k^T,G_k^T)\ $ if
$n=2k+1$\\
with $\varepsilon=0$&
$(I_k,J_k(-1))\oplus
(I_k,J_k(-1))\ $ if
$n=2k$ $(k$ is
odd$)$\\\hline
$\begin{matrix}
\text{Matrix
\eqref{cm2}}\\
\text{with
$\varepsilon=1$}
\end{matrix}$
&
$(I_n,J_n((-1)^{n+1})$
\\\hline
Matrix \eqref{cm3}&
$(I_{2k},J_{2k}(1))\oplus
(I_{2k},J_{2k}(1))$\\\hline
\end{tabular}
\end{equation}
\end{theorem}

Let $\mathbb F$ be an
algebraically closed
field with nonidentity
involution. Fix
$i\in\mathbb F$ such
that $i^2=-1$. It is
known (see Lemma
\ref{l_real}) that
each element of
$\mathbb F$ is
uniquely representable
in the form $a+bi$
with $a,b$ in $
\mathbb P:=\{\lambda
\in\mathbb F\,|\,\bar
\lambda=\lambda\}$,
and the involution on
$\mathbb F$ is
``complex
conjugation'':
$\overline{a+bi}=a-bi$.
Moreover, $\mathbb P$
is ordered and
$a^2+b^2$ has a unique
positive real root,
which is called the
\emph{modulus} of
$a+bi$ and is denoted
by $|a+bi|$.

\begin{theorem}
\label{t_matr1}
{\rm(a)} Each square
matrix $A$ over an
algebraically closed
field $\mathbb F$ with
nonidentity involution
is *congruent to a
direct sum, determined
uniquely up to
permutation of
summands, of
tridiagonal matrices
of two types:
\begin{equation}\label{cmi1}
\begin{bmatrix}
    0&1 &&&0\\
     \lambda&0&1 \\
    &\lambda &0&\ddots&\\
    &&\ddots&\ddots&1 \\
    0&&&\lambda &0
  \end{bmatrix}_{n},\quad
  \begin{matrix}
\lambda\in\mathbb F,\
 |\lambda| \ne 1,\\
    \text{each nonzero
    $\lambda$ is
determined}\\
\text{ up to
replacement by
$\bar\lambda^{-1}$,}\\
\text{$\lambda= 0$ if
$n$
is odd}\\
  \end{matrix}
\end{equation}
$($one can take
$|\lambda|<1$ if $n$
is even$)$; and
\begin{equation}\label{cmi2}
\mu
\begin{bmatrix}
    1 &1&&&&0\\
    -1&0&1&\\
    &1&0&1&\\
  &&-1&0&1&\\
&&&1&0&\ddots&\\
  0&&&&  \ddots&\ddots
  \end{bmatrix}_{n},
    \qquad
    \mu\in\mathbb F,\
  |\mu |=1.
\end{equation}

{\rm(b)} The Kronecker
canonical form of
$(A^*,A)$ under
equivalence determines
the direct sum
asserted in {\rm(a)}
uniquely up to
permutation of
summands and
multiplication of any
direct summand of type
\eqref{cmi2} by $-1$.
For each direct
summand $B$ of type
\eqref{cmi1} or
\eqref{cmi2}, the
Kronecker canonical
form of $(B^*,B)$ is
given in the following
table:
\begin{equation}\label{tab2}
\renewcommand{\arraystretch}{1.2}
\begin{tabular}{|c|c|}
\hline $B$& Kronecker
canonical form of $
(B^*,B)$\\
\hline\hline Matrix
\eqref{cmi1}&
$\begin{array}{rl}
(F_k,G_k)\oplus
(F_k^T,G_k^T)
&\text{if $n=2k+1$}\\
(J_k(\bar\lambda),I_k)
\oplus
(I_k,J_k(\lambda))
&\text{if $n=2k$}
\end{array}$
\\\hline
Matrix \eqref{cmi2}&
$(I_n,J_n((-1)^{n+1}
\bar\mu^{-1}\mu)
)$\\\hline
\end{tabular}
\end{equation}
\end{theorem}

\section{Four lemmas}
\label{s1a}

In this section we
prove four lemmas that
we use in later
sections. In the first
lemma we collect known
results about
algebraically closed
fields with
involution; i.e., a
bijection $a\mapsto
\bar{a}$ satisfying
$\overline{a+b}=\bar
a+ \bar b$,
$\overline{ab}=\bar a
\bar b$ and $\bar{\bar
a}=a$.

\begin{lemma}\label{l_real}
Let\/ $\mathbb F$ be
an algebraically
closed field with
nonidentity involution
$\lambda\mapsto\bar\lambda$,
and let
\begin{equation}\label{123}
\mathbb
P:=\bigl\{\lambda
\in{\mathbb
F}\,\bigr|\,
\bar{\lambda }=\lambda
\bigr\}.
\end{equation}
Then $\mathbb F$ has
characteristic $0$,
\begin{equation}\label{1pp11}
\mathbb F={\mathbb
P}+{\mathbb P}i,\qquad
i^2=-1,
\end{equation}
and the involution has
the form
\begin{equation}\label{1ii}
\overline{a+bi}=a-bi,\qquad
a,b\in\mathbb P.
\end{equation}
Moreover, the field\/
${\mathbb P}$ has a
unique linear ordering
$\le$ such that
\begin{equation*}\label{slr}
\text{$a>0$ and\,
$b>0$}
 \quad\Longrightarrow\quad
\text{$a+b>0$ and\,
$ab>0$}.
\end{equation*}
The positive elements
of\/ $\mathbb P$ with
respect to this
ordering are the
squares of nonzero
elements.  Every
algebraically closed
field of
characteristic $0$
possesses a
nonidentity
involution.
\end{lemma}

\begin{proof}
If $\mathbb F$ is an
algebraically closed
field with nonidentity
involution $\lambda
\mapsto \bar{\lambda
}$, then this
involution is an
automorphism of order
2. Hence ${\mathbb F}$
has degree $2$ over
the field ${\mathbb
P}$ defined in
\eqref{123}. By
Corollary 2 in
\cite[Chapter VIII, \S
9]{len}, $\mathbb P$
has characteristic $0$
and every element of
${\mathbb F}$ is
uniquely representable
in the form $a+bi$
with $a,b\in{\mathbb
P}$. Since the
involution is an
automorphism of
${\mathbb F}$,
$\bar{i}^2=-1$. So
$\bar{i}=-i$ and the
involution is
\eqref{1ii}. Due to
Proposition 3 in
\cite[Chapter XI, \S
2]{len}, $\mathbb P$
is a real closed
field, and so the
statements about the
ordering $\le$ follow
from \cite[Chapter XI,
\S 2, Theorem 1]{len}.
By \cite[\S 82,
Theorem 7c]{wan},
every algebraically
closed field of
characteristic $0$
contains at least one
real closed subfield
and hence it can be
represented in the
form \eqref{1pp11} and
possesses the
involution
\eqref{1ii}.
\end{proof}

The canonical form
problem for pairs of
symmetric or
skew-symmetric
matrices under
congruence reduces to
the canonical form
problem for matrix
pairs under
equivalence due to the
following lemma, which
was proved in
\cite[\S\,95, Theorem
3]{mal} for complex
matrices. Roiter
\cite{roi} (see also
\cite{ser1dir,ser_izv})
extended this lemma to
arbitrary systems of
linear mappings and
bilinear forms over an
algebraically closed
field of
characteristic
different from $2$.

\begin{lemma}
\label{l_mal} Let $(A,
B)$ and $(A',B')$ be
given pairs of
$n\times n$ matrices
over an algebraically
closed field $\mathbb
F$ of characteristic
different from $2$.
Suppose that $A$ and
$A'$ are either both
symmetric or both
skew-symmetric, and
also that $B$ and $B'$
are either both
symmetric or both
skew-symmetric. Then
$(A, B)$ and $(A',B')$
are congruent if and
only if they are
equivalent.
\end{lemma}

\begin{proof}
If $(A, B)$ and
$(A',B')$ are
congruent then they
are equivalent.

Conversely, let $(A,
B)$ and $(A',B')$ be
equivalent; i.e.,
$R^TAS=A'$ and
$R^TBS=B'$ for some
nonsingular $R$ and
$S$. Then
\[
R^TAS=A' =\varepsilon
(A')^T=\varepsilon
S^TA^TR=S^TAR,
\]
in which $\varepsilon
=1$ if $A$ and $A'$
are symmetric and
$\varepsilon =-1$ if
$A$ and $A'$ are
skew-symmetric. Write
$M:=SR^{-1}$. Then
\[
AM=M^TA,\quad
AM^2=(M^T)^2A,\ \ldots
\]
and so $Af(M)=f(M)^TA$
for every polynomial
$f\in\mathbb F[x]$. If
there exists
$f\in\mathbb F[x]$
such that $f(M)^2=M$,
then for $N:=f(M)R$ we
have
\[
N^TAN=R^Tf(M)^TAf(M)R
=R^TAf(M)^2R=R^TAMR
=R^TAS=A'.
\]
Repeating the argument
for the matrix $B$, we
obtain $N^TBN=B'$.
Consequently, $(A, B)$
and $(A',B')$ are
congruent.

It remains to find
$f\in\mathbb F[x]$
such that $f(M)^2=M$.
 Let
\[
(x-\lambda_1)^{k_1}\cdots
(x-\lambda_t)^{k_t},\qquad
\lambda_i\ne\lambda_j
\text{ if }i\ne j,
\]
be the characteristic
polynomial of $M$. We
can reduce $M$ to
Jordan canonical form
and obtain
\begin{equation*}\label{kmh}
M=J_1\oplus\dots
\oplus J_t,\qquad
J_i=\lambda_iI_{k_i}
+F_i,\quad
F_i^{k_i}=0.
\end{equation*}
For the polynomial
\[\varphi_i(x):
=\prod_{j\ne
i}(x-\lambda
_j)^{k_j}\] we have
\begin{equation}\label{hhg}
\varphi_i(M)
=0_{k_1+\dots+k_{i-1}}
\oplus\varphi_i(J_i)
\oplus
0_{k_{i+1}+\dots+k_t}
\end{equation}
($0_k$ denotes the
$k\times k$ zero
matrix). The field
$\mathbb F$ is
algebraically closed
of characteristic not
$2$, all $\lambda_i$
and $\varphi_i(\lambda
_i)$ are nonzero, so
for each $i=1,\dots,t$
there exist
polynomials $\psi_i,
\tau_i\in\mathbb F[x]$
such that
\begin{equation*}\label{kjh}
\psi_i(x)^2\equiv
\lambda _i+x,\quad
\varphi_i(\lambda
_i+x)\tau_i(x)\equiv
\psi_i(x) \mod x^{k_i}
\end{equation*}
(the coefficients of
$\psi_i$ and $\tau_i$
are determined
successively from
these congruences).
Then
$f(x):=\sum_i\varphi
_i(x)\tau_i(x-\lambda
_i)$ is the required
polynomial. Indeed, by
\eqref{hhg}
\[
f(M)=
\bigoplus_i\varphi
_i(J_i)\tau_i(J_i-\lambda
_iI_{k_i})=\bigoplus_i\varphi
_i(\lambda
_iI_{k_i}+F_i)\tau_i(F_i)=
\bigoplus_i\psi_i(F_i)
\]
and so
$$
f(M)^2=\bigoplus_i
\psi_i(F_i)^2=
\bigoplus_i (\lambda_i
I_{k_i}+F_i)=
\bigoplus_i J_i=M.
\eqno\qedhere
$$
\end{proof}

For each matrix of the
form
\begin{equation}
\label{vdv}
A=\begin{bmatrix}
\varepsilon &a_1&&&&0\\
a'_1&0&b_1\\
&b'_1&0&a_2\\
&&a'_2&0&b_2\\
&&&b'_2&0&\ddots\\
0&&&&\ddots&\ddots
  \end{bmatrix}_{n},
\end{equation}
          define
\begin{equation*}
\label{vdv1} {\cal
P}(A):=\begin{bmatrix}
b_k &a'_{k}&&&&&0\\
&\ddots&\ddots\\
&&b_1&a'_1\\
&&&\varepsilon&a_1 \\
&&&&b'_1&\ddots\\
&&&&&\ddots&a_{k} \\
0&&&&&&b'_k
  \end{bmatrix}_{2k+1}
  \quad
\text{if $n=2k+1$}
\end{equation*}
       and
\begin{equation*}
\label{vdv2} {\cal
P}(A):=\begin{bmatrix}
a_k &b'_{k-1}&&&&&&0\\
&a_{k-1}&\ddots\\
&&\ddots&b'_1\\
&&&a_1&\varepsilon \\
&&&&a'_1&b_1\\
&&&&&a'_2&\ddots\\
&&&&&&\ddots&b_{k-1} \\
0&&&&&&&a'_k
  \end{bmatrix}_{2k}\quad
\text{if $n=2k$.}
\end{equation*}

\begin{lemma}
\label{l_per} Every
pair $(A,B)$ of
$n\times n$ matrices
of the form
\eqref{vdv} is
equivalent to $({\cal
P}(A),{\cal P}(B))$.
\end{lemma}

\begin{proof}
If $n=2k+1$, then we
rearrange rows
$1,2,\dots,2k+1$ in
$A$ and in $B$ as
follows:
\[
2k,\ 2k-2,\ \dots,\
2,\ 1,\ 3,\ \dots,\
2k-1,\ 2k+1,
\]
and their columns in
the inverse order:
\[
2k+1,\ 2k-1,\ \dots,\
3,\ 1,\ 2,\ \dots,\
2k-2,\ 2k.
\]
If $n=2k$, then we
rearrange the rows of
$A$ and $B$ as
follows:
\[
2k-1,\ 2k-3,\ \dots,\
3,\ 1,\ 2,\ \dots,\
2k-2,\ 2k,
\]
and their columns in
the inverse order:
\[
2k,\ 2k-2,\ \dots,\
2,\ 1,\ 3,\ \dots,\
2k-3,\ 2k-1.
\]
The pair that we
obtain is $({\cal
P}(A),{\cal P}(B))$.
\end{proof}

For a sign $\sigma
\in\{+,-\}$ and a
nonnegative integer
$k$, define the
$2k$-by-$2k$ matrix
\begin{equation*}\label{ndyd}
M^{\sigma }_k
:=\begin{bmatrix}
  0 & 1 \\
  \sigma
 1 & 0
\end{bmatrix}\oplus\dots
\oplus \begin{bmatrix}
  0 & 1\\
  \sigma
 1 & 0
\end{bmatrix}\qquad
\text{($k$ summands).}
\end{equation*}
Thus, $M^{\sigma }_0$
is $0$-by-$0$.

\begin{lemma}
\label{l_ge} Let
$\sigma ,\tau
\in\{+,-\}$ and
$k\in\mathbb N$. Then
the following pairs
are equivalent:
\begin{align}
 \label{2}
(0_1\oplus M^{\sigma
}_{k},\: M^{\tau}
_{k}\oplus
0_1)&\approx
(F_{k},G_{k})
 \oplus
(F_{k}^T,G_{k}^T),
     \\
\label{1} (I_1\oplus
M^{\sigma }_{k},\:
M^{\tau} _{k}\oplus
0_1)&\approx
(I_{2k+1},J_{2k+1}(0)),
    \\ \label{4}
(0_1\oplus
M^{\sigma}_{k-1}\oplus
0_1,\:
M^{\tau}_{k})&\approx
(J_k(0),I_k)\oplus
(J_k(0),I_k),
    \\ \label{3}
(I_1\oplus M^{\sigma
}_{k-1}\oplus 0_1,\:
M^{\tau}_{k})&\approx
(J_{2k}(0),I_{2k}).
\end{align}
\end{lemma}

\begin{proof}
Let $\varepsilon
\in\{0,1\}.$ By Lemma
\ref{l_per},
\begin{equation*}\label{jdg}
([\varepsilon]\oplus
M^{\sigma }_{k},
M^{\tau }_{k}\oplus
0_1)\approx (I_{k}
\oplus[\varepsilon]
\oplus I_{k},\
J_{2k+1}(0)),
\end{equation*}
which proves \eqref{2}
and \eqref{1}, and
\[
([\varepsilon]\oplus
M^{\sigma
}_{k-1}\oplus 0_1,
M^{\tau}_{k})\approx
\left(\begin{bmatrix}
    0&1&&&&0\\
    &0&\cdot&&&\\
    &&\cdot&\varepsilon &&\\
    &&&\cdot&\cdot&
    \\
    &&&&\cdot&1\\
    0&&&&&0
  \end{bmatrix},\
  \begin{bmatrix}
    1&&&&&0\\
    &1&&&&\\
    &&\cdot&&&\\
    &&&\cdot&&
    \\
    &&&&\cdot&\\
    0&&&&&1
  \end{bmatrix}\right),
\] which proves
\eqref{4} and
\eqref{3}.
\end{proof}

\section{Pairs of
symmetric matrices}
\label{s1}

In this section, we
give two tridiagonal
canonical forms of
pairs of symmetric
matrices under
congruence.

\subsection{First
canonical form}
\label{sub1}

\begin{theorem}\label{th1a}
{\rm(a)} Over an
algebraically closed
field $\mathbb F$ of
characteristic
different from $2$,
every pair $(A,B)$ of
symmetric matrices of
the same size is
congruent to a direct
sum, determined
uniquely up to
permutation of
summands, of
tridiagonal pairs of
two types:
    %%%%ss1n
\begin{equation}\label{sss1n}
\left(\begin{bmatrix}
    0 &1&&&&0\\
    1&0&0\\
    &0&0&1\\
    &&1&0&0\\
    &&&0&0&\ddots\\
    0&&&&\ddots&\ddots\\
  \end{bmatrix}_{n},\
  \begin{bmatrix}
    \varepsilon &
    \lambda &&&&0\\
    \lambda &0&1\\
    &1&0&\lambda \\
    &&\lambda &0&1\\
    &&&1 &0&\ddots\\
    0&&&&\ddots&\ddots
  \end{bmatrix}_{n}\right),
\ \
 \begin{matrix}
\lambda\in\mathbb F,\\
   \varepsilon
\in\{0,1\},
  \end{matrix}
\end{equation}
in which
$\varepsilon=1$ if $n$
is even and
$\lambda=0$ if $n$ is
odd; and
      %%%%ss2n
\begin{equation}\label{ss2n}
\left(\begin{bmatrix}
    1&0&&&&0\\
    0&0 &1&&&\\
    &1&0&0\\
       &&0&0&1\\
 &&&1&0&\ddots\\
    0&&&&\ddots&\ddots \\
  \end{bmatrix}_{n},\
  \begin{bmatrix}
    \lambda &1&&&&0\\
    1&0&\lambda &&&\\
    &\lambda &0&1\\
        &&1 &0&\lambda\\
    &&&\lambda&0&\ddots\\
    0&&&&\ddots&\ddots
  \end{bmatrix}_{n}\right),
\end{equation}
in which $\lambda=0$
if $n$ is even and
$\lambda \in\mathbb F$
if $n$ is odd.

{\rm(b)} This direct
sum is determined
uniquely up to
permutation of
summands by the
Kronecker canonical
form of $(A,B)$ under
equivalence. The
Kronecker canonical
form of each of the
direct summands is
given in the following
table:
\begin{equation}\label{tab3}
\renewcommand{\arraystretch}{1.2}
\begin{tabular}{|c|c|}
\hline
 Pair& Kronecker
canonical form of
the pair
\\
\hline\hline
\eqref{sss1n}&
$\begin{array}{rl}
(F_k,G_k)
 \oplus
(F_k^T,G_k^T)
&\text{if $n=2k+1$ and
$\varepsilon=0$}\\
(J_{n}(0),I_{n})
&\text{if $n$ is odd
and
$\varepsilon=1$}\\
(I_{n},J_{n}(\lambda))
&\text{if $n$ is even}
\end{array}$
\\\hline
\eqref{ss2n}&
$\begin{array}{rl}
(I_{n},
J_{n}(\lambda))
&\text{if $n$ is odd}
\\
(J_{n}(0),I_{n})
&\text{if $n$ is even}
\end{array}$\\\hline
\end{tabular}
\end{equation}
\end{theorem}

\begin{proof}
Let the Kronecker
canonical form of
$(A,B)$ be
\begin{equation*}\label{ksc}
\bigoplus_i
(I_{m_i},J_{m_i}
(\lambda_i))
  \oplus\bigoplus_j
(J_{n_j}(0),I_{n_j})
\oplus
 \bigoplus_l
 (F_{s_l},G_{s_l})
 \oplus
 \bigoplus_r
(F_{t_r}^T,G_{t_r}^T).
\end{equation*}
Since $A$ and $B$ are
symmetric,
\[
(A,B)\approx\bigoplus_i
(I_{m_i},J_{m_i}
(\lambda_i))
  \oplus\bigoplus_j
(J_{n_j}(0),I_{n_j})
\oplus
 \bigoplus_l
 (F_{s_l}^T,G_{s_l}^T)
 \oplus
 \bigoplus_r
(F_{t_r},G_{t_r}).
\]
Thus, we can make
$s_1=t_1,\
s_2=t_2,\dots$ by
reindexing $\{t_r\}$,
and obtain that the
Kronecker canonical
form of $(A,B)$ is
\begin{equation}\label{kscw}
\bigoplus_i
(I_{m_i},J_{m_i}
(\lambda_i))
  \oplus\bigoplus_j
(J_{n_j}(0),I_{n_j})
\oplus
 \bigoplus_l
 \Big((F_{s_l},G_{s_l})
 \oplus
(F_{s_l}^T,G_{s_l}^T)\Big).
\end{equation}
This sum is determined
by $(A,B)$ uniquely up
to permutation of
summands. In view of
Lemma \ref{l_mal}, it
remains to prove
\eqref{tab3}.

The pair \eqref{sss1n}
with $n=2k+1$ and
$\varepsilon=0$ has
the form
$(M^+_{k}\oplus
0_1,0_1\oplus
M^+_{k})$; by
\eqref{2} it is
equivalent to
$(F_k,G_k)
 \oplus
(F_k^T,G_k^T)$.

The pair \eqref{sss1n}
with $n=2k+1$ and
$\varepsilon=1$ has
the form $
(M^+_{k}\oplus
0_1,\:I_1\oplus
M^+_{k})$; by
\eqref{1} it is
equivalent to
$(J_{2k+1}(0),I_{2k+1})$.

The pair \eqref{sss1n}
with $n=2k$ has the
form $(M^+_k,\:\lambda
M^+_k+(I_1\oplus
M^+_{k-1}\oplus
0_1))$; it is
equivalent to
$(I_{2k}, \lambda
I_{2k}+ J_{2k}(0))
=(I_{2k},J_{2k}(\lambda))$
since \eqref{3}
ensures that
\begin{equation}\label{jdd1}
(M^+_k,\: I_1\oplus
M^+_{k-1}\oplus
0_1)\approx
(I_{2k},J_{2k}(0)).
\end{equation}

The pair \eqref{ss2n}
with $n=2k+1$ has the
form $(I_1\oplus
M^+_{k},\:\lambda
(I_1\oplus M^+_{k})+
(M^+_{k}\oplus 0_1))$;
by \eqref{1} it is
equivalent to
$(I_{2k+1},
J_{2k+1}(\lambda))$.

The pair \eqref{ss2n}
with $n=2k$ has the
form $(I_1\oplus
M^+_{k-1}\oplus
0_1,\:M^+_k)$; by
\eqref{3} it is
equivalent to
$(J_{2k}(0),I_{2k})$.
\end{proof}

\subsection{Second
canonical form}
\label{sub2}

In this section, we
give another
tridiagonal canonical
form of pairs of
symmetric matrices for
congruence. This form
is not a direct sum of
tridiagonal matrices
of the form
\eqref{kdu}. It is
based on the
\raisebox{1.5pt}{-}\!\!Dokovi\'{c}
and Zhao's tridiagonal
canonical form of
symmetric matrices for
orthogonal similarity
\cite{dok-zhao_tridiag}
and resembles the
Kronecker canonical
form of matrix pairs
for equivalence.

For each positive
integer $n$, let $N_n$
denote any fixed
$n\times n$
tridiagonal symmetric
matrix over $\mathbb
F$ that is similar to
$J_n(0)$. Following
\cite[p.\,79]{dok-zhao_tridiag},
we can take as $N_n$
the value
$N(a_1,\dots,a_n,b)$
of the polynomial
matrix
\begin{equation*}\label{kgk}
N(x_1,\dots,x_n,y):=\begin{bmatrix}
 x_1&y&&0\\
 y&x_2&\ddots&\\
 &\ddots&\ddots&y\\
  0&&y&x_n
  \end{bmatrix}
\end{equation*}
at any nonzero
solution
$(a_1,\dots,a_n,b)\in
\mathbb F^{n+1}$ of
the system
\[
c_1(x_1,\dots,x_n,y)=0,\
\dots,\
c_n(x_1,\dots,x_n,y)=0
\]
of equations whose
left parts are the
coefficients of the
characteristic
polynomial
$t^n+c_1t^{n-1}+
\dots+c_n$ of
$N(x_1,\dots,x_n,y)$.
Then $0$ is the only
eigenvalue of $N_n$,
$b\ne 0$, $\rank
N_n=n-1$, and $N_n$ is
similar to $J_n(0)$.

If $\mathbb F$ has the
characteristic $0$,
then
\cite[p.\,81]{dok-zhao_tridiag}
ensures that we can
also take
\begin{equation*}\label{luf}
N_n=\begin{bmatrix}
 n-1&id_1&&&0\\
 id_1&n-3&id_2&&\\
 &id_2&n-5&\ddots&\\
  &&\ddots&\ddots&id_{n-1}\\
   0&&&id_{n-1}&1-n
\end{bmatrix},\quad
\begin{matrix}
d_l:=\sqrt{l(n-l)},\\
i^2=-1.
\end{matrix}
\end{equation*}

\begin{theorem}\label{th1b}
Over an algebraically
closed field $\mathbb
F$ of characteristic
different from $2$,
every pair $(A,B)$ of
symmetric matrices of
the same size is
congruent to a direct
sum, determined
uniquely up to
permutation of
summands, of
tridiagonal pairs of
three types:
\begin{equation}\label{lyg}
(I_n,\lambda I_n+N_n)
\text{ with $\lambda
\in\mathbb F$};\qquad
(N_n, I_n);
\end{equation}
and
\begin{equation}\label{ssnew}
\left(\begin{bmatrix}
    0&1&&&&0\\
    1&0 &0&&&\\
    &0&0&1\\
       &&1&0&0\\
 &&&0&0&\ddots\\
    0&&&&\ddots&\ddots \\
  \end{bmatrix}_{2k+1},\
  \begin{bmatrix}
    0 &0&&&&0\\
    0&0&1 &&&\\
    &1 &0&0\\
        &&0 &0&1\\
    &&&1&0&\ddots\\
    0&&&&\ddots&\ddots
  \end{bmatrix}_{2k+1}\right).
\end{equation}

{\rm(b)} This direct
sum is determined
uniquely up to
permutation of
summands by the
Kronecker canonical
form of $(A,B)$ for
equivalence. The
Kronecker canonical
form of each of the
direct summands is
given in the following
table:
\begin{equation}\label{tab4}
\renewcommand{\arraystretch}{1.2}
\begin{tabular}{|c|c|}
\hline
 Pair& Kronecker
canonical form of the
pair
\\
\hline\hline
 $(I_n,\lambda
I_n+N_n)$&$(I_{n},J_{n}
(\lambda))$\\\hline
$(N_n,
I_n)$&$(J_{n}(0),I_{n})$
\\\hline
\eqref{ssnew}& $
(F_{k},G_{k})
 \oplus
(F_{k}^T,G_{k}^T)$\\\hline
\end{tabular}
\end{equation}
\end{theorem}

\begin{proof}
In view of
\eqref{kscw} and Lemma
\ref{l_mal}, it
suffices to prove
\eqref{tab4}. The
equivalences
\[
(I_{n},\lambda
I_{n}+N_{n}) \approx
(I_{n},J_{n}
(\lambda))\ \text{ and
}\ (N_{n}, I_{n})
\approx
(J_{n}(0),I_{n})
\]
are valid since $N_n$
is similar to
$J_n(0)$. The pair
\eqref{ssnew} is
\eqref{sss1n} with
$n=2k+1$ and
$\varepsilon=0$; by
\eqref{tab3} it is
equivalent to $
(F_{k},G_{k})
 \oplus
(F_{k}^T,G_{k}^T)$.
\end{proof}

\section{Pairs of
matrices, in which the
first is symmetric and
the second is
skew-symmetric}
\label{s2}

\begin{theorem}\label{t_sc}
Over an algebraically
closed field $\mathbb
F$ of characteristic
different from $2$,
every pair $(A,B)$ of
matrices of the same
size, in which $A$ is
symmetric and $B$ is
skew-symmetric, is
congruent to a direct
sum, determined
uniquely up to
permutation of
summands, of
tridiagonal pairs of
three types:
\begin{equation}\label{sc1}
  \left(\begin{bmatrix}
    0&1&&&0\\
    1&0&1\\
    &1&0&\ddots
    \\
    &&\ddots&\ddots&1\\
    0&&&1&0
  \end{bmatrix}_{2k},\
  \begin{bmatrix}
    0&\lambda&&&0\\
    -\lambda&0&\lambda\\
    &-\lambda&0&\ddots
    \\
    &&\ddots&\ddots&\lambda\\
    0&&&-\lambda&0
  \end{bmatrix}_{2k}\right),
  \quad
  \begin{matrix}
\lambda\in\mathbb F,\\
\lambda \ne 0,
  \end{matrix}
\end{equation}
in which $\lambda$ is
determined up to
replacement by
$-\lambda$;

\begin{equation}\label{sc2}
\left(\begin{bmatrix}
    \varepsilon &0&&&&0\\
    0&0&1\\
    &1&0&0\\
    &&0&0&1\\
    &&&1&0&\ddots\\
    0&&&&\ddots&\ddots
  \end{bmatrix}_n,\
  \begin{bmatrix}
    0&1&&&&0\\
    -1&0&0\\
    &0&0&1\\
    &&-1&0&0\\
    &&&0&0&\ddots\\
    0&&&&\ddots&\ddots
  \end{bmatrix}_n\right),
\end{equation}
in which $\varepsilon
=1$ if $n$ is a
multiple of $4$, and
$\varepsilon
\in\{0,1\}$ otherwise;
and
\begin{equation}\label{sc3}
\left(
\begin{bmatrix}
    0 &1&&&&0\\
    1&0&0\\
    &0&0&1\\
    &&1&0&0\\
    &&&0&0&\ddots\\
    0&&&&\ddots&\ddots\\
  \end{bmatrix}_{4k},\
   \begin{bmatrix}
    0&0&&&&0\\
    0&0&1\\
    &-1&0&0\\
    &&0&0&1\\
    &&&-1&0&\ddots\\
   0&&&&\ddots&\ddots
\end{bmatrix}_{4k}\right).
\end{equation}

{\rm(b)} This direct
sum is determined
uniquely up to
permutation of
summands by the
Kronecker canonical
form of $(A,B)$ under
equivalence. The
Kronecker canonical
form of each of the
direct summands is
given in the following
table:
\begin{equation}\label{tab5}
\renewcommand{\arraystretch}{1.2}
\begin{tabular}{|c|c|}
\hline
 Pair& Kronecker
canonical form of the
pair
\\
\hline\hline
\eqref{sc1}&$(I_k,J_k(\lambda
))\oplus
(I_k,J_k(-\lambda))$
with $\lambda \ne 0$
   \\\hline
\eqref{sc2} with
$\varepsilon=0$
     &
$\begin{array}{rl}
(F_k,G_k)
 \oplus
(F_k^T,G_k^T)
&\text{if
$n=2k+1$}\\
(J_k(0),I_k)\oplus
(J_k(0),I_k) &\text{if
$n=2k$  $(k$ is
odd$)$}
\end{array}$\\\hline
\eqref{sc2} with
$\varepsilon=1$
     &
$\begin{array}{rl}
(I_n,J_n(0)) &\text{if
$n$ is odd}\\
(J_n(0),I_n) &\text{if
$n$ is even}
\end{array}$\\\hline
\eqref{sc3}&
$(I_{2k},J_{2k}(0))\oplus
(I_{2k},J_{2k}(0))$\\\hline
\end{tabular}
\end{equation}
\end{theorem}

\begin{proof}
The Kronecker
canonical form of
$(A,B)$ is a direct
sum of pairs of the
types:
\begin{itemize}
  \item[\rm(i)]
$(I_k,J_k(\lambda
))\oplus
(I_k,J_k(-\lambda))$,
in which $\lambda\ne
0$ if $k$ is odd,

  \item[\rm(ii)]
$(I_n,J_n(0))$ with
odd $n$,

  \item[\rm(iii)]
$(J_k(0),I_k)\oplus
(J_k(0),I_k)$  with
odd $k$,

  \item[\rm(iv)]
$(J_n(0),I_n)$  with
even $n$,

  \item[\rm(v)]
$(F_k,G_k)
 \oplus
(F_k^T,G_k^T)$.
\end{itemize}
This statement was
proved in
\cite[Section 4]{thom}
for pairs of complex
matrices and goes back
to Kronecker's 1874
paper; see the
historical remark at
the end of Section 4
in \cite{thom}. The
proof remains valid
for matrix pairs over
$\mathbb F$ (or see
\cite[Theorem
4]{ser_izv}).

In view of Lemma
\ref{l_mal}, it
suffices to prove
\eqref{tab5}.

By Lemma \ref{l_per},
\eqref{sc1} is
equivalent to
\[
(J_k(1),-\lambda
J_k(-1))\oplus
(J_k(1),\lambda
J_k(-1)),
\]
which is equivalent to
(i) with $\lambda\ne
0$.

The pair \eqref{sc2}
with $n=2k+1$ has the
form $([\varepsilon]
\oplus M^+_{k},
M^-_{k}\oplus 0_1)$;
by \eqref{2} and
\eqref{1} this pair is
equivalent to (v) if
$\varepsilon=0$ or
(ii) if
$\varepsilon=1$.

The pair \eqref{sc2}
with $n=2k$ has the
form
\begin{equation}\label{hdt}
([\varepsilon]\oplus
M^+_{k-1}\oplus 0_1,
M^-_{k}),
\end{equation}
in which $\varepsilon
\in\{0,1\}$ if $k$ is
odd and $\varepsilon
=1$ if $k$ is even.
Due to \eqref{4} and
\eqref{3}, \eqref{hdt}
is equivalent to (iii)
if $\varepsilon=0$ or
to (iv) if
$\varepsilon=1$.

The pair \eqref{sc3}
has the form
$(M^+_{2k}, 0_1\oplus
M^-_{2k-1}\oplus
0_1)$, and by
\eqref{4} it is
equivalent to (i) with
$\lambda =0$ and $k$
replaced by $2k$.
\end{proof}

\section{Pairs of
skew-symmetric
matrices} \label{ss1}

\begin{theorem}\label{t_skew}
Over an algebraically
closed field $\mathbb
F$ of characteristic
different from $2$,
every pair $(A,B)$ of
skew-symmetric
matrices of the same
size is congruent to a
direct sum, determined
uniquely up to
permutation of
summands, of
tridiagonal pairs of
two types:
\begin{equation}\label{cc1}
\left(
\begin{bmatrix}
    0&1&&&&0\\
    -1&0&0\\
    &0&0&1\\
    &&-1&0&0\\
    &&&0&0&\ddots\\
    0&&&&\ddots&\ddots
  \end{bmatrix}_{2k},\
  \begin{bmatrix}
    0&\lambda  &&&&0\\
    -\lambda  &0&1\\
    &-1&0&\lambda\\
    &&-\lambda&0&1\\
    &&&-1  &0&\ddots\\
    0&&&&\ddots&\ddots
  \end{bmatrix}_{2k}\right)
\end{equation}
and
   %%cc23
\begin{equation}\label{cc23}
\left(\begin{bmatrix}
    0 &0&&&&0\\
    0&0&1\\
    &-1&0&0\\
    &&0&0&1\\
    &&&-1&0&\ddots\\
    0&&&&\ddots&\ddots
  \end{bmatrix}_n,\
  \begin{bmatrix}
    0&1&&&&0\\
    -1&0&0\\
    &0&0&1\\
    &&-1&0&0\\
    &&&0&0&\ddots\\
    0&&&&\ddots&\ddots
  \end{bmatrix}_n\right)
\end{equation}
in which
$k,n\in\mathbb N$ and
$\lambda \in\mathbb
F$.

{\rm(b)} This direct
sum is determined
uniquely up to
permutation of
summands by the
Kronecker canonical
form of $(A,B)$ under
equivalence. The
Kronecker canonical
form of each of the
direct summands is
given in the following
table:
\begin{equation}\label{tab6}
\renewcommand{\arraystretch}{1.2}
\begin{tabular}{|c|c|}
\hline
 Pair& Kronecker
canonical form of the
pair
\\
\hline\hline
\eqref{cc1}&$(I_k,J_k(\lambda
))\oplus
(I_k,J_k(\lambda))$
   \\\hline
\eqref{cc23}
     &
$\begin{array}{rl}
(F_k,G_k)
 \oplus
(F_k^T,G_k^T)
&\text{if
$n=2k+1$}\\
(J_k(0),I_k)\oplus
(J_k(0),I_k) &\text{if
$n=2k$}
\end{array}$\\\hline
\end{tabular}
\end{equation}
\end{theorem}

\begin{proof}
The Kronecker
canonical form of
$(A,B)$ under
equivalence is a
direct sum of pairs of
three types:
\begin{gather*}\label{ksqw}
((I_k,J_k
(\lambda))\oplus
(I_k,J_k
(\lambda))),\qquad
((J_k(0),I_k)\oplus
(J_k(0),I_k)),
    \\
 ((F_k,G_k)
 \oplus
(F_k^T,G_k^T)).
\end{gather*}
This statement was
proved in
\cite[Section 4]{thom}
for pairs of complex
matrices, but the
proof remains valid
for pairs over
$\mathbb F$ (or see
\cite[Theorem
4]{ser_izv}). In view
of Lemma \ref{l_mal},
it suffices to prove
\eqref{tab6}.

The pair \eqref{cc1}
has the form
$(M^-_k,\lambda
M^-_k+(0_1\oplus
M^-_{k-1}\oplus 0_1))$
and by \eqref{4} it is
equivalent to
\[
(I_{k},\lambda I_k+
J_{k} (0))\oplus
(I_{k},\lambda I_k+
J_{k} (0)) =
(I_{k},J_{k} (\lambda
))\oplus (I_{k},J_{k}
(\lambda )).
\]

The pair \eqref{cc23}
with $n=2k+1$ has the
form $ (0_1\oplus
M^-_{k},M^-_{k}\oplus
0_1)$; by \eqref{2} it
is equivalent to
$(F_{k},G_{k})
 \oplus
(F_{k}^T,G_{k}^T)$.

The pair \eqref{cc23}
with $n=2k$ has the
form $(0_1\oplus
M^-_{k-1}\oplus
0_1,M^-_k)$; by
\eqref{4} it is
equivalent to
$(J_{k}(0),I_{k})\oplus
(J_{k}(0),I_{k})$.
\end{proof}

\section{Matrices with
respect to congruence}
\label{s_pr}

In this section we
prove Theorem
\ref{t_matr}.
\medskip

(a) Each square matrix
$A$ can be expressed
uniquely as the sum of
a symmetric and a
skew-symmetric matrix:
\begin{equation*}\label{pair16}
A=A_{\text{sym}}
+A_{\text{sk}},\qquad
A_{\text{sym}}:
=\frac{A+A^{T}}2
,\quad
A_{\text{sk}}:=\frac{A-A^{T}}2.
\end{equation*}
Two matrices $A$ and
$B$ are congruent if
and only if the
corresponding pairs
$(A_{\text{sym}},
A_{\text{sk}})$ and
$(B_{\text{sym}},
B_{\text{sk}})$ are
congruent. Therefore,
adding the first and
the second matrices in
each of the canonical
pairs from Theorem
\ref{t_sc} gives three
types of canonical
matrices for
congruence:
\begin{equation}\label{cm1a}
  \begin{bmatrix}
    0&1+\mu &&0\\
    1-\mu &\ddots&\ddots&\\
    &\ddots&0&1
    +\mu \\
    0&&1-\mu &0
  \end{bmatrix}_{2k},\qquad
  \begin{matrix}
    \text{$\mu\ne 0,$}\\
    \text{$\mu$ is
determined up}\\
\text{to replacement
by $-\mu$};
  \end{matrix}
\end{equation}
\eqref{cm2}; and
\eqref{cm3}. We can
assume that $\mu \ne
-1$ because the
congruence
transformation
\begin{equation}\label{hfe}
X\mapsto S^TXS,\qquad
S:=\begin{bmatrix}0&&1\\
&\ddd&\\1&&0
\end{bmatrix},
\end{equation}
maps \eqref{cm1a} with
$\mu= -1$ into
\eqref{cm1a} with
$\mu= 1$. If we
multiply all the odd
columns and rows of
\eqref{cm1a} by
$(1+\mu)^{-1}$ (this
is a transformation of
congruence), we obtain
\eqref{cm1} with
\begin{equation}\label{hvf}
\lambda=\frac{1-\mu}
{1+\mu}.
\end{equation}

The parameter $\mu$ is
determined up to
replacement by $-\mu$,
so each $\lambda\ne 0$
is determined up to
replacement by
$\lambda^{-1}$,
whereas $\lambda=0$ is
determined uniquely
since it corresponds
to $\mu = 1$ and we
assume that $\mu \ne
-1$. We have
$\lambda\ne \pm 1$
because $\mu \ne 0$
and ${-1+\mu}\ne
{1+\mu}$. The
parameter $\lambda$ is
an arbitrary element
of $\mathbb F$ except
for $\pm 1$ since
substituting
$\mu=(1-\lambda)/
(1+\lambda)$ into
\eqref{hvf} gives the
identity.
\medskip

(b) Let $A$ be the
matrix \eqref{cm1}. By
Lemma \ref{l_per}, the
pair $(A^T,A)$ is
equivalent to
\begin{equation}\label{hhy}
\left(
\left[\begin{array}
{c|c}
\begin{matrix}
\lambda &1\\&\lambda
&\ddots\\&&\ddots&1\\
&&&\lambda
\end{matrix}
&0\\ \hline 0&
  \begin{matrix}
1 &\lambda\\&1
&\ddots\\
&&\ddots&\lambda\\
&&&1
\end{matrix}
\end{array}
\right]\!,
\left[\begin{array}
{c|c}
\begin{matrix}
1 &\lambda\\&1
&\ddots\\
&&\ddots&\lambda\\
&&&1\end{matrix} &0\\
\hline 0&
  \begin{matrix}
\lambda &1\\&\lambda
&\ddots\\&&\ddots&1\\
&&&\lambda
\end{matrix}
\end{array}
\right] \right)\!,
\end{equation}
which is equivalent to
$(J_k(\lambda), I_k)
\oplus
(I_k,J_k(\lambda))$
since $\lambda \ne \pm
1$. This verifies the
assertion about the
matrix \eqref{cm1} in
table \eqref{tab1}.

The remaining
assertions about the
matrices \eqref{cm2}
and \eqref{cm3} in
table \eqref{tab1}
follow from the
corresponding
assertions about the
matrices \eqref{sc2}
and \eqref{sc3} in
table \eqref{tab5}:
the matrices
\eqref{cm2} and
\eqref{cm3} have the
form $A=B+C$ in which
$(B,C)$ is \eqref{sc2}
or \eqref{sc3}, and so
$(A^T,A)=(B-C,B+C)$.
For example, if $A$ is
\eqref{cm2} with
$\varepsilon =1$, then
by \eqref{tab5}
\[
(B,C)\approx
  \begin{cases}
    (I_n,J_n(0))
    & \text{if $n$ is odd}, \\
        (J_n(0),I_n)
& \text{if $n$ is
even},
  \end{cases}
\]
and we have
\[
(A^T,A)\approx
  \begin{cases}
  (I_n-J_n(0), I_n+J_n(0))
  \approx
    (I_n,J_n(1))
    & \text{if $n$ is odd},
       \\
(J_n(0)-I_n, J_n(0)+
I_n)
  \approx
(I_n,J_n(-1)) &
\text{if $n$ is even}.
  \end{cases}
\]
The proof of Theorem
\ref{t_matr} is
complete.

\section{Matrices with
respect to
*congruence}
\label{s_pr1}

In this section we
prove Theorem
\ref{t_matr1}.

Let $\mathbb F$ be an
algebraically closed
field with nonidentity
involution represented
in the form
\eqref{1pp11}. A
canonical form of a
square matrix $A$ over
$\mathbb F$ for
*congruence was given
in \cite{ser_izv} and
was improved in
\cite{hor-ser_transp}
(a direct proof that
the matrices in
\cite{hor-ser_transp}
are canonical is given
in
\cite{hor-ser_regul,
hor-ser2}): $A$ is
*congruent to a direct
sum, determined
uniquely up to
permutation of
summands, of matrices
of three types:
\begin{equation}\label{eqq}
\begin{bmatrix}0&I_k\\
J_k(\lambda) &0
\end{bmatrix}\
(\lambda\ne 0,\
|\lambda |\ne 1),\quad
\mu\begin{bmatrix}
0&&&1
\\
&&\ddd&i\\
&1&\ddd&\\
1&i&&0
\end{bmatrix}\ (|\mu|=1),
 \quad
J_n(0),
\end{equation}
in which $\lambda$ is
determined up to
replacement by
$\bar\lambda^{-1}$. It
follows from the proof
of Theorem 3 in
\cite{ser_izv} that
instead of \eqref{eqq}
one can take any set
of matrices
\begin{equation*}\label{vdf}
P_{2k}(\lambda),\qquad
\mu Q_n,\qquad J_n(0)
\end{equation*}
(with the same
conditions on $\lambda
$ and $\mu $) such
that
\begin{equation}\label{azs1}
(P_{2k}(\lambda)^*,
P_{2k}(\lambda))
\approx
(J_k(\bar\lambda),
I_k) \oplus
(I_k,J_k(\lambda))
\end{equation}
and
\begin{equation}\label{azs2}
(Q_n^*, Q_n) \approx
(I_n,J_n(\nu_n)),
\end{equation}
in which
$\nu_1,\nu_2,\dots$
are any elements of
$\mathbb F$ with
modulus one.

\begin{proof}[Proof of
Theorem \ref{t_matr1}]

Let $P_{2k}(\lambda)$
be the matrix
\eqref{cmi1} with
$\lambda \ne 0$ and
let $Q_n$ be the
matrix \eqref{cmi2}
with $\mu=1$. Since
the matrix
\eqref{cmi1} with
$\lambda=0$ is $
J_n(0)$, it suffices
to prove that
\eqref{azs1} and
\eqref{azs2} are
fulfilled.

By Lemma \ref{l_per},
$(P_{2k}(\lambda)^*,
P_{2k}(\lambda))$ is
equivalent to the pair
\eqref{hhy} with
$\bar\lambda$ instead
of $\lambda$ in the
first matrix. This
proves \eqref{azs1}
since $|\lambda |\ne
1$.

The matrix $Q_n$ is
\eqref{cm2} with
$\varepsilon =1$. Due
to \eqref{tab1},
\[
(Q_n^*,Q_n)=(Q_n^T,Q_n)
\approx
(I_n,J_n((-1)^{n+1});
\]
this ensures
\eqref{azs2} with
$\nu_n:=(-1)^{n+1}$.

The assertion about
the matrix
\eqref{cmi1} with
$\lambda=0$ in table
\eqref{tab2} follows
from the equivalence
\[
(J_n(0)^T,J_n(0))\approx
  \begin{cases}
(F_k,G_k)\oplus
(F_k^T,G_k^T)
&\text{if $n=2k+1$,}\\
(J_k(0),I_k) \oplus
(I_k,J_k(0)) &\text{if
$n=2k$},
  \end{cases}
\]
which was established
in the proof of
Theorem 3 in
\cite{ser_izv}.
\end{proof}

\section{Pairs of
Hermitian matrices}
\label{s_her}

\begin{theorem}\label{t_her}
{\rm(a)} Over an
algebraically closed
field $\mathbb F$ with
nonidentity involution
represented in the
form \eqref{1pp11},
every pair $(A,B)$ of
Hermitian matrices of
the same size is
*congruent to a direct
sum, determined
uniquely up to
permutation of
summands, of
tridiagonal pairs of
two types:
     %%%he1
\begin{equation}\label{he1}
\left(\begin{bmatrix}
    0&1&&&0\\
    1&0&1\\
    &1&0&\ddots\\
    &&\ddots&\ddots&1\\
    0&&&1&0
    \end{bmatrix}_{n},\
  \begin{bmatrix}
    0&\mu &&&0\\
    \bar\mu &0&\mu \\
    &\bar\mu &0&\ddots\\
    &&\ddots&\ddots&\mu \\
    0&&&\bar\mu &0
    \end{bmatrix}_{n}\right),
\end{equation}
in which
$\mu\in\mathbb
F\smallsetminus\mathbb
P$ if $n$ is even,
$\mu=\pm i$ if $n$ is
odd, and $\mu$ is
determined up to
replacement by
$\bar\mu$; and
   %%he2
\begin{equation}\label{he2}
\left(\begin{bmatrix}
    a &b&&&&0\\
    b&0&a\\
    &a&0&b\\
    &&b&0&a\\
    &&&a&0&\ddots\\
    0&&&&\ddots&\ddots
  \end{bmatrix}_n,\
  \begin{bmatrix}
    b&-a&&&&0\\
    -a&0&b\\
    &b&0&-a\\
    &&-a&0&b\\
    &&&b&0&\ddots\\
    0&&&&\ddots&\ddots
  \end{bmatrix}_n\;\right),
\end{equation}
in which
$a,b\in\mathbb P$ and
$a^2+b^2=1$.

{\rm(b)} The Kronecker
canonical form of
$(A,B)$ under
equivalence determines
this direct sum
uniquely up to
permutation of
summands and
multiplication by $-1$
any direct summand of
type \eqref{he2}. The
Kronecker canonical
form of each of the
direct summands is
given in the following
table:
\begin{equation}\label{tab2a}
\renewcommand{\arraystretch}{1.2}
\begin{tabular}{|c|c|}
\hline
 Pair& Kronecker
canonical form of the
pair
\\
\hline\hline
\eqref{he1}&
$\begin{array}{rl}
(F_k,G_k)\oplus
(F_k^T,G_k^T)
&\text{if $n=2k+1$}\\
(I_k,J_k(\mu)) \oplus
(I_k,J_k(\bar\mu))
&\text{if $n=2k$}
\end{array}$
\\\hline
\eqref{he2}&
$\begin{array}{rl}
(I_n,J_n( b/a))
&\text{if $n$ is odd
and $a\ne 0$}\\
(I_n,J_n(-a/b))
&\text{if $n$ is even
and $b\ne 0$}\\
(J_n(0),I_n)
&\text{otherwise}
\end{array}$
\\\hline
\end{tabular}
\end{equation}
\end{theorem}

\begin{proof}
(a) Each square matrix
$A$ over $\mathbb F$
has a
\textit{Cartesian
decomposition}
\begin{equation*}
\label{pair6}
 A=B+iC,
\qquad
B:=\frac{A+A^*}{2}
,\quad C:=\frac{i(
A^*-A)}{2},
\end{equation*}
in which both $B$ and
$C$ are Hermitian. Two
square matrices $A$
and $A'$ are
*congruent if and only
if the corresponding
pairs $(B,C)$ and
$(B',C')$ are
*congruent. Therefore,
if we apply the
Cartesian
decomposition to the
canonical matrices for
*congruence from
Theorem \ref{t_matr1},
we obtain canonical
pairs of Hermitian
matrices for
*congruence. To
simplify these
canonical pairs, we
multiply \eqref{cmi1}
by $2$ (this is a
transformation of
*congruence), and
using \eqref{1pp11}
take $\mu$ in
\eqref{cmi2} to have
the form $a+bi$ with
$a,b\in\mathbb P$.
Thus, every pair
$(A,B)$ of Hermitian
matrices of the same
size is *congruent to
a direct sum,
determined uniquely up
to permutation of
summands, of pairs of
two types:
    %%%hemm1
\begin{equation}\label{hemm1}
\left(\begin{bmatrix}
    0&\!1+\bar\lambda
    \!&&0\\
    \!1+\lambda\!&0
    &\ddots\\
    &\ddots&\ddots&\!1
    +\bar\lambda\! \\
    0&&\!1+\lambda\! &0
  \end{bmatrix}_{n},\
  i\begin{bmatrix}
    0 &
    \!\bar\lambda-1\!&&0\\
   \! 1-\lambda\!
   &0&\ddots\\
    &\ddots&\ddots&
    \!\bar\lambda-1\! \\
    0&&\!1-\lambda\! &0
  \end{bmatrix}_{n}\right),
\end{equation}
in which
$\lambda\in\mathbb F$,
$|\lambda| \ne 1$,
each nonzero $\lambda$
is determined up to
replacement by
$\bar\lambda^{-1}$,
and $\lambda= 0$ if
$n$ is odd; and
   %%%hemm2
\begin{equation}
    \label{hemm2}
\left(\begin{bmatrix}
    a&bi&&&&0\\
 -bi&0&a\\
   &a&0&bi\\
   &&-bi&0&a\\
 &&&a&0&\ddots\\
    0&&&&\ddots&\ddots
  \end{bmatrix}_n,\
  \begin{bmatrix}
    b&-ai&&&&0\\
 ai&0&b\\
   &b&0&-ai\\
   &&ai&0&b\\
 &&&b&0&\ddots\\
    0&&&&\ddots&\ddots
  \end{bmatrix}_n\;\right),
\end{equation}
in which $a^2+b^2=1$.

Let us prove that the
pairs \eqref{hemm1}
and \eqref{hemm2} are
*congruent to the
pairs \eqref{he1} and
\eqref{he2}.

We obtain \eqref{he2}
if we apply the
*congruence
transformation
$X\mapsto S^*XS$ with
\[
S:=\diag(1,-i,-i,-1,-1,
i,i,1,1,-i,-i,-1,-1,\ldots)
\]
to the matrices of
\eqref{hemm2}.

The pair \eqref{hemm1}
with $\lambda=0$  is
the pair \eqref{he1}
with $\mu=-i$, which
is *congruent to
\eqref{he1} with
$\mu=i$ via the
transformation
\eqref{hfe}.

It remains to consider
\eqref{hemm1} with
$\lambda\ne 0$. Then
$n$ is even. Applying
to the matrices of
\eqref{hemm1} the
*congruence
transformation
$X\mapsto S^*XS$ with
\[
S:=\diag\left(1,\
\frac{1}{1+\bar\lambda},\
\frac{1+\lambda}{1
+\bar\lambda},\
\frac{1+\lambda}{(1
+\bar\lambda)^2},\
\frac{(1+\lambda)^2}{(1
+\bar\lambda)^2},\
\frac{(1+\lambda)^2}{(1
+\bar\lambda)^3},\
\ldots \right),
\]
(the denominator is
nonzero since
$|\lambda| \ne 1$), we
obtain \eqref{he1}
with
\begin{equation}\label{nkl}
\mu:=\frac{\bar\lambda-1}{
\bar\lambda+1}i.
\end{equation}
Since $\lambda$ is
nonzero and is
determined up to
replacement by
$\bar\lambda ^{-1}$,
we have that $\mu\ne
-i$ and $\mu$ is
determined up to
replacement by
\[
\frac{\lambda^{-1}-1}
{\lambda^{-1}+1}i=
\frac{1-\lambda}
{1+\lambda}i=\bar\mu.
\]
Every $\mu\in\mathbb
F$ except for $i$ can
be represented in the
form \eqref{nkl} with
$\lambda=(i-\bar\mu)/
(i+\bar\mu)$. We do
not impose the
condition $\mu\ne\pm
i$ in \eqref{he1}
because \eqref{he1}
with $\mu=\pm i$ is
*congruent to
\eqref{hemm1} with
$\lambda=0$.

Let us prove that the
condition
$|\lambda|\ne 1$ is
equivalent to the
condition $\mu\notin
\mathbb P$. If
$|\lambda|= 1$ and
$\lambda=a+bi \ne-1$
with $a,b\in\mathbb
P$, then
\begin{equation}\label{yyy}
\mu=\frac{(\bar\lambda-1)
(\lambda+1)}{(\bar\lambda
+1)(\lambda+1)}i=
\frac{\bar\lambda\lambda-
\lambda+ \bar\lambda
-1}
{\bar\lambda\lambda
+\lambda+\bar\lambda
+1}i= \frac{-bi}{1+a}i
\in\mathbb P.
\end{equation}
Each $\mu\in\mathbb P$
can be represented in
the form \eqref{yyy}
as follows:
$\mu=b/(1+a)$, in
which
\[
a:=\frac{1-\mu^2}{1+\mu^2}
\quad\text{and}\quad
b:=\frac{2\mu}{1+\mu^2}
\qquad (\text{then
}a^2+b^2=1).
\]

(b) Lemma \ref{l_per}
ensures the assertion
about the pair
\eqref{he1} in table
\eqref{tab2a}.

The pair \eqref{he2}
has the form $(aX+bY,
bX-aY)$, in which
$(X,Y)$ is
\eqref{ss2n} with
$\lambda=0$. By
\eqref{tab3},
$(X,Y)\approx
(I_n,J_n(0))$ if $n$
is odd, and
$(X,Y)\approx
(J_n(0),I_n) $ if $n$
is even. Therefore,
\[
\text{Pair
\eqref{he2}}\approx
  \begin{cases}
(aI_n+bJ_n(0),bI_n-aJ_n(0))
 & \text{if $n$
is odd}, \\
(aJ_n(0)+bI_n,bJ_n(0)-aI_n)
& \text{if $n$ is
even}.
  \end{cases}
\]
This validates the
assertion about the
pair \eqref{he2} in
table \eqref{tab2a}.
\end{proof}

\begin{remark}\label{rem2}
The pair \eqref{he2}
with two dependent
parameters, which was
obtained from the
Cartesian
decomposition of
\eqref{cmi2}, can be
replaced by $0$- and
$1$-parameter matrices
as follows. The
matrices \eqref{cmi2}
have the form $\mu A$,
in which $\mu=a+bi$,
$a,b\in\mathbb P$, and
$a^2+b^2=1$. If
$\mu\ne\pm i$, then
$a\ne 0$ and $\mu A$
is *congruent to
$|a|^{-1}\mu
A=\pm(1+ci)A$ with
$c\in\mathbb P$. Now
apply the Cartesian
decomposition to $\pm
iA$ and $ \pm(1+ci)A$
with $c\in\mathbb P$.
\end{remark}

\end{document}